\documentclass[letterpaper,10pt,conference]{ieeeconf}%

\usepackage{hyperref} 
\usepackage{graphicx}
\usepackage{amsthm}
\usepackage{algorithm}
\usepackage[noend]{algpseudocode}
\usepackage[cmex10]{amsmath}
\usepackage{amsfonts,amssymb}
\usepackage{times}
\usepackage{bm}
\usepackage{bbm}
\usepackage{color}
\usepackage{graphicx}
\usepackage{subcaption}
\usepackage{amsfonts}
\usepackage{cite}
\usepackage{mathtools}
\usepackage{amssymb}%
\usepackage{tabu}
\usepackage{color}
\providecommand{\U}[1]{\protect\rule{.1in}{.1in}}

\IEEEoverridecommandlockouts
\newtheorem{theorem}{Theorem}

\newtheorem{lemma}{Lemma}

\newtheorem{definition}{Definition}
\theoremstyle{definition}
\newtheorem{remark}{Remark}
\newtheorem{assumption}{Assumption}

\usepackage{tcolorbox}
\usepackage{lipsum}
\tcbuselibrary{skins,breakable}
\usetikzlibrary{shadings,shadows}
\usepackage{xcolor}

    {\endtcolorbox}

\makeatletter
\def\BState{\State\hskip-\ALG@thistlm}
\makeatother

\begin{document}

\title{{\LARGE \textbf{Optimal Threshold-Based Control Policies for Persistent
Monitoring on Graphs}}}
\author{Nan Zhou$^{1}$, Christos G. Cassandras$^{1,2}$, Xi Yu$^{3}$, and Sean B.
Andersson$^{1,3}$\\{\small $^{1}$Division of Systems Engineering, $^{2}$Department of Electrical
and Computer Engineering, $^{3}$Department of Mechanical Engineering}\\Boston University, Boston, MA 02215, USA\\E-mail:\texttt{\{nanzhou,cgc,xyu,sanderss\}@bu.edu 
\thanks{* The work of
Cassandras and Zhou is supported in part by NSF under grants ECCS-1509084,
CNS-1645681, and IIP-1430145, by AFOSR under grant FA9550-15-1-0471, by DOE
under grant DOE-46100, by MathWorks and by Bosch. The work of Andersson and Yu
is supported in part by NSF through grants ECCS-1509084 and CMMI-1562031.} }}
\maketitle

\begin{abstract}
We consider the optimal multi-agent persistent monitoring problem defined by a
team of cooperating agents visiting a set of nodes (targets) on a graph with
the objective of minimizing a measure of overall node state uncertainty. The
solution to this problem involves agent trajectories defined both by the
sequence of nodes to be visited by each agent and the amount of time spent at
each node. Since such optimal trajectories are generally intractable, we
propose a class of distributed threshold-based parametric controllers through
which agent transitions from one node to the next are controlled by threshold
parameters on the node uncertainty states. The resulting behavior of the
agent-target system can be described by a hybrid dynamic system. This enables
the use of Infinitesimal Perturbation Analysis (IPA) to determine on line
(locally) optimal threshold parameters through gradient descent methods and
thus obtain optimal controllers within this family of threshold-based
policies. We further show that in a single-agent case the IPA gradient is
monotonic, which implies a simple structure whereby an agent visiting a node
should reduce the uncertainty state to zero before moving to the next node.
Simulation examples are included to illustrate our results and compare them to
optimal solutions derived through dynamic programming when this is possible.

\end{abstract}


\section{Introduction}

\label{sec:intro}

The cooperative multi-agent \emph{persistent monitoring}\ problem arises when
agents are tasked to monitor a dynamically changing environment which cannot
be fully covered by a stationary agent allocation. Thus, persistent monitoring
differs from traditional consensus \cite{ren2005survey} and coverage control
\cite{zhong2011distributed} problems due to the continuous need to explore
changes in the environment. In many cases, this exploration process leads to the
discovery of various \textquotedblleft points of interest\textquotedblright,
which, once detected, become \textquotedblleft data sources\textquotedblright%
\ or \textquotedblleft targets\textquotedblright\ that need to be perpetually
monitored. This paradigm applies to surveillance systems, such as when a team
must monitor large regions for changes, intrusions, or other dynamic events
\cite{michael2011persistent}, or when it is responsible for sampling and
monitoring environmental parameters such as temperature
\cite{leonard2010coordinated}. It also finds use in particle tracking in
molecular biology where the goal is to track multiple individual biological
macromolecules to understand their dynamics and their interactions
\cite{ashley2016tracking,cromer2011tracking}. In contrast to sweep coverage and
patrolling \cite{Smith:2012fq,lin2015optimal}, the problem we address here
focuses on a \emph{finite number} of data sources or \textquotedblleft
targets\textquotedblright\ (typically larger than the number of agents). The
goal of the agent team is to collect information from each target so as to
reduce a metric of uncertainty about its state. This uncertainty naturally
increases while no agent is present in its vicinity and decreases when it is
being monitored (or \textquotedblleft sensed\textquotedblright) by one or more
agents. Thus, the objective is to minimize an overall measure of target
uncertainty by controlling the movement of all agents.

Our previous work \cite{zhou2016optimal} considered the persistent monitoring
problem in a one-dimensional space, formulated it as an optimal control
problem and showed that the solution can be reduced to a parametric controller
form. In particular, the optimal agent trajectories are characterized by a
finite number of points where each agent switches direction and by a dwell
time at each such point. However, in two-dimensional (2D) spaces, it has been
shown that such parametric representations for optimal agent trajectories no
longer hold \cite{lin2015optimal}. Nonetheless, various forms of parametric
trajectories (e.g., ellipses, Lissajous curves, interconnected linear
segments) can still be near-optimal or at least offer an alternative
\cite{lin2015optimal,zhou2017graphPM}. These approaches limit agent
trajectories to certain forms which, while they possess desirable properties (e.g.,
periodicity), cannot always capture the dynamic changes in target
uncertainties and may lead to poor local optima
\cite{lahijanian2010motion,cassandras2013optimal,zhou2016optimal}.

In this paper, we take a different direction in the 2D persistent monitoring
problem. Rather than parameterizing agent trajectories, we adopt a more
abstract point of view whereby targets are nodes in a graph and their
connectivity defines feasible agent trajectories along the edges of the graph.
A trajectory is specified by a sequence of nodes and an associated dwell time
at each node in the sequence. In this setting, there is a travel time defined
for each edge in the graph which is determined in advance according to the
actual target topology and has the added benefit of accounting for constraints
such as physical obstacles in the 2D space which the graph is designed to
avoid. The controller associated with each agent determines $(i)$ the dwell
time at the current node and $(ii)$ the next node to be visited with the goal
of optimizing a given performance metric. The complexity in this optimization
problem is significant \cite{yu2017optimal} and the presence of real-valued
dwell time decision variables makes it much harder than that of Traveling
Salesman Problems \cite{bektas2006multiple}, which are already computationally
intensive and do not scale well. Since deriving such optimal trajectories is
generally intractable, we consider a class of controllers based on a set of
\emph{threshold parameters} {associated with} the target uncertainties, hence
taking into account the time-varying nature of target states. Thus, the
parameterization in this graph-based setting is imposed on the {target}
thresholds rather than on the shape of agent trajectories in the underlying 2D
Euclidean environment. By adjusting the thresholds, we can control the agent
behavior in terms of target visiting and dwelling and, therefore, optimize a
given performance metric within the specific parametric controller family
considered. From a modeling standpoint, this results in a hybrid dynamic
system whose state consists of agent positions and target uncertainties. From
an optimization standpoint, the goal is to determine optimal thresholds
(parameters) that minimize a given metric.

The contribution of the paper lies in the graph-based setup of the 2D
persistent monitoring problem, the formulation of a threshold-based parametric
optimization problem, and a solution approach based on Infinitesimal
Perturbation Analysis (IPA) \cite{cassandras2010perturbation} to determine on
line the gradient of the objective function and to obtain (possibly local)
optimal threshold parameters through gradient descent. As we will see,
optimizing these thresholds not only affects the dwell time that the agent
should spend at each node, but also naturally adjusts and seeks to optimize
the node visiting sequence. Our approach is \emph{distributed} since the
decisions made by an agent at some node are based on uncertainty states of
neighboring nodes only. Moreover, we exploit the event-driven nature of IPA to
also render it \emph{scalable} in the number of \emph{events} in the system
and not the state space (in contrast to solutions dependent on dynamic
programming). An additional contribution is to show that in the case of a
one-agent system the IPA gradient is monotonic in the thresholds involved
which implies a simple optimal structure: the agent visiting a node should
reduce the uncertainty state to zero before moving to the next node. This is
consistent with a similar earlier result established in \cite{yu2017optimal}.

The paper is organized as follows. Section II formulates the 2D persistent
monitoring problem on a graph and introduces the parametric family of
threshold-based agent controllers we subsequently analyze. Section III
provides a solution of the optimization problem obtained through event-driven
IPA gradient estimation. In Section IV, we present our analysis of the
one-agent case with the key result that all optimal dwell times are specified
by zero threshold values. Section V includes simulation examples, including
comparisons with optimal solutions derived through dynamic programming when
this is possible. Section VI concludes the paper.

\section{Problem formulation}

\label{sec:PM_form} Consider $N$ agents and $M$ targets in a 2D mission space.
The agent positions are $s_{a}(t)\in\mathbb{R}^{2}$, $a=1,\dots,N$ and the
target locations are $X_{i}\in\mathbb{R}^{2}$, $i=1,\dots,M$.

\textbf{Target uncertainty model}. Following the model in
\cite{zhou2016optimal}, we define uncertainty functions $R_{i}(t)$ associated
with targets $i=1,\ldots,M$, with the following properties: $(i)$ $R_{i}(t)$
increases with a prespecified rate $A_{i}$ if no agent is visiting it, $(ii)$
$R_{i}(t)$ decreases with a rate $B_{i}N_{i}(t)$ where $B_{i}$ is the rate at
which an agent collects data from target $i$, hence decreasing its uncertainty
state, and $N_{i}(t)=\sum_{a=1}^{N}\mathbf{1}\{s_{a}(t)=X_{i}\}$ is the number
of agents at target $i$ at time $t$, and $(iii)$ $R_{i}(t)\geq0$ for all $t$.
We model the target uncertainty state dynamics as follows:
\begin{equation}
\hspace*{-2mm}\dot{R}_{i}(t)=\left\{
\begin{array}
[c]{ll}%
0 &\text{if }R_{i}(t)=0\text{ and }A_{i}\leq B_{i}N_{i}(t)  \\
A_{i}-B_{i}N_{i}(t)& \text{otherwise } 
\end{array}
\right.  \label{eq:DynR}%
\end{equation}
This model has an attractive queueing system interpretation as explained in
\cite{zhou2016optimal}, where each target is associated with an
\textquotedblleft uncertainty queue\textquotedblright\ with input rate $A_{i}$
and service rate $B_{i}N_{i}(t)$ controllable through the agent movement. Note
that compared with the model in \cite{zhou2016optimal}, where each agent has a
finite sensing range $r_{a}$ allowing it to decrease $R_{i}(t)$ as long as
$\left\Vert s_{a}(t)-X_{i}\right\Vert <r_{a}$, here the agent's sensing range
is ignored and the joint detection probability of a target by agents is
replaced by the summation $N_{i}(t)$ above. This is done for simplicity to
accommodate the graph topology we will adopt; the analysis can be extended to
the original model in \cite{zhou2016optimal} at the expense of added notation
and the use of a sensing model for each agent.

\textbf{Agent model}. The position of agent $a$ is denoted by $s_{a}%
(t)=[x_{a}(t),y_{a}(t)]^{\top}$ for $a=1,\ldots,N$ and its dynamics in 2D are
given by:
\begin{equation}
\dot{s}_{a}(t)=[v_{a}(t)\cos(u_{a}(t)),v_{a}(t)\sin(u_{a}(t))]^{\top}
\label{eq:DynA}%
\end{equation}
where the agent's velocity is scaled and bounded such that $\Vert
v_{a}(t)\Vert\leq1$ and the agent's heading is $u_{a}(t)\in\lbrack0,2\pi)$.

\textbf{Objective function}. Our goal is to determine the optimal control
(both $v_{a}(t)$ and $u_{a}(t)$) for all agents under which the average
uncertainty metric in \eqref{eq:costfunction} across all targets is minimized
over a given time horizon $T$. Setting $\mathbf{v}(t)=[v_{1}(t),\ldots
,v_{N}(t)]$ and $\mathbf{u}(t)=[u_{1}(t),\ldots,u_{N}(t)]$, we aim to solve
the following optimal control problem:
\begin{equation}
\mathbf{P1:}\min_{\mathbf{v}(t),\mathbf{u}(t)}\text{ \ }J=\frac{1}{T}\int%
_{0}^{T}\sum_{i=1}^{M}R_{i}(t)\,dt \label{eq:costfunction}%
\end{equation}
subject to target dynamics in \eqref{eq:DynR} and agent dynamics in \eqref{eq:DynA}.

Obtaining a complete solution of \textbf{P1} generally requires solving a
computationally hard Two Point Boundary Value Problem (TPBVP) which amounts to
a 2D functional search in both $v_{a}(t)$ and $u_{a}(t)$ for each agent over
$t\in\lbrack0,T]$. Unlike the 1D case in \cite{zhou2016optimal}, the problem
cannot be reduced to a parametric one as shown in \cite{lin2015optimal}.
However, it is still easy to show that the optimal agent speed is limited to
$\Vert v_{a}^{\star}(t)\Vert\in\{1,0\}$ depending on whether the agent is
dwelling or traveling. In particular, the Hamiltonian associated with
\textbf{P1} is
\begin{equation}%
\begin{split}
H=  &  \sum_{i=1}^{M}R_{i}(t)+\sum_{i=1}^{M}\lambda_{i}(t)\dot{R}_{i}(t)\\
&  +\sum_{a=1}^{N}v_{a}\left(  \lambda_{a}^{x}(t)\cos u_{a}(t)+\lambda_{a}%
^{y}(t)\sin u_{a}(t)\right)
\end{split}
\end{equation}
and a straightforward application of the Pontryagin minimum principle implies
that $v_{a}^{\ast}(t)= \pm 1$ depending on the sign of
$\left(  \lambda_{a}^{x}(t)\cos u_{a}(t)+\lambda_{a}^{y}(t)\sin u_{a}%
(t)\right)  $, or $v_{a}^{\ast}(t)=0$ in singular arcs that may exist. The
analysis is similar to the one in \cite{lin2015optimal} and is, therefore,
omitted here.

Using this optimal control structure and the underlying target topology, we
make a further simplification by constraining agent movements to a graph
$G=(V,E)$ where the set of vertices (nodes) is defined by an indexed list of
targets $V=\{1,\ldots,M\}$ and the set of edges (links) $E$ contains all
feasible direct connections between them. Note that if there are obstacles in
the underlying space, we can introduce \textquotedblleft way
points\textquotedblright\ to define feasible paths avoiding the obstacles,
where a way point $j$ is included in the set $V$ with an associated
uncertainty state $R_{j}(t)=0$ for all $t\geq0$. In this graph topology, the
agent headings $u_{a}(t)$ are limited to the finite set $V$, i.e.,
$u_{a}(t)\in V=\{1,\ldots,M\}$.

Therefore, we have reduced \textbf{P1} to a simpler problem of determining
$(i)$ the dwell time for each agent at each node when $v_{a}(t)=0$ and $(ii)$
the control (heading) $u_{a}(t)$ when $v_{a}(t)\neq0$. The complete state of
this system is defined by $\mathbf{s}(t)=[s_{1}(t),\ldots,s_{N}(t)]$ and
$\mathbf{R}(t)=[R_{1}(t),\ldots,R_{M}(t)]$ so that the control should be
expressed as $u_{a}(\mathbf{s}(t),\mathbf{R}(t))$.

\begin{figure}[pt]
\centering
\includegraphics[width=0.95\linewidth]{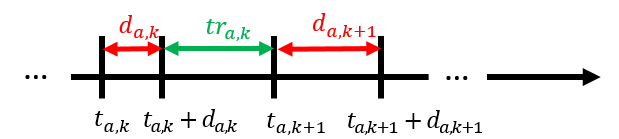}\caption{{\protect\small An
agent control trajectory: $d_{a,k}$ is the $k$-th dwell time and $tr_{a,k}$ is
the $k$-th travel time.}}%
\label{fig:Control_traj}%
\end{figure}

Figure \ref{fig:Control_traj} shows a typical control trajectory and helps
pinpoint the behavior of each agent controller. The trajectory consists of a
sequence of intervals $[t_{a,k},t_{a,k+1})$ where the agent's node visits are
indexed by $k=1,2,\ldots$ and $t_{a,k}$ is the time of the $k$-th visit at any
node. This interval contains $t_{a,k}+d_{a,k}$, the time when the agent leaves
the current node it is visiting. Note that on an optimal trajectory:%
\[
t_{a,k+1}^{\ast}=t_{a,k}^{\ast}+d_{a,k}^{\ast}+\Vert s_{a}^{\ast}%
(t_{a,k})-s_{a}^{\ast}(t_{a,k+1})\Vert
\]
since the optimal agent speed when transitioning between nodes satisfies
$\Vert v_{a}^{\star}(t)\Vert=1$. Thus, the agent controller's role when
visiting some node $i$ is to determine the optimal dwelling time
$d_{a,k}^{\ast}$ and next node $u_{a}^{\ast}(t_{a,k}^{\ast}+d_{a,k}^{\ast})$.
Clearly, $u_{a}(t)$ is a piecewise constant right-continuous function of time
and the values of $u_{a}(t)$ belong to the set
\begin{equation}
u_{a}\left(  t\right)  \in\{i\}\cup\mathcal{N}_{i}\quad\text{ if }%
s_{a}(t)=X_{i} \label{eq:control_dwell}%
\end{equation}
where $\mathcal{N}_{i}$ is the neighborhood of node $i$ defined as follows.

\begin{definition}
\label{def:neighboring_targets} The neighborhood of node $i$ is the set
$\mathcal{N}_{i}=\{j:(i,j)\in E,$ $j\in V\}$.
\end{definition}

Since the control $u_{a}\left(  t\right)  $ switches only at times
$t_{a,k}+d_{a,k}$ (see Fig. \ref{fig:Control_traj}), let us concentrate on a
time interval $[t_{a,k},t_{a,k}+d_{a,k})$ during which $s_{a}(t)=X_{i}$ for
some $i\in V$. Observe that for $t\geq t_{a,k}$ either $u_{a}\left(  t\right)
=i$ or it switches to a new value $j\in\mathcal{N}_{i}$. The condition under
which such a switch occurs may generally be expressed as $g_{i,j}%
(\mathbf{s}(t),\mathbf{R}(t))\leq0$, i.e., $g_{i,j}(\mathbf{s}(t),\mathbf{R}%
(t))$ is a switching function associated with a transition from node $i$ to
node $j$. Let us define%
\[
\tau_{a,k}^{j}=\inf_{t\geq t_{a,k}}\{g_{i,j}(\mathbf{s}(t),\mathbf{R}(t))=0\}
\]
and set%
\[
t_{a,k}+d_{a,k}=\min_{j\in\mathcal{N}_{i}}\{\tau_{a,k}^{j}\}
\]
so that the change in the agent's node assignment occurs at the earliest time
that one of the switching functions satisfies $g_{i,j}(\mathbf{s}%
(t),\mathbf{R}(t))=0$. Thus, the task of the controller is to determine
optimal switching functions $g_{i,j}^{\ast}(\mathbf{s}(t),\mathbf{R}(t))$ for
all $j\in\mathcal{N}_{i}$ whenever $s_{a}(t)=X_{i}$ and then evaluate
$\min_{j\in\mathcal{N}_{i}}\{\tau_{a,k}^{j}\}$ to specify the optimal dwelling
time $d_{a,k}^{\ast}$. Therefore,
\begin{align}
u_{a}^{\ast}\left(  t\right)   &  =i\text{, \, }t\in\lbrack t_{a,k}^{\ast
},t_{a,k}^{\ast}+d_{a,k}^{\ast})\label{u_optimal}\\
u_{a}^{\ast}\left(  t_{a,k}^{\ast}+d_{a,k}^{\ast}\right)   &  =\arg\min
_{j\in\mathcal{N}_{i}}\{\tau_{a,k}^{j}\}\nonumber
\end{align}
In effect, whenever $s_{a}(t)=X_{i}$, the state space defined by all feasible
values of $[\mathbf{s}(t),\mathbf{R}(t)]$ is partitioned into $\left\vert
\mathcal{N}_{i}\right\vert +1$ regions, denoted by $\mathcal{R}_{i}$ and
$\mathcal{R}_{j}$, $j\in\mathcal{N}_{i}$. The controller keeps the agent at
node $i$ as long as $[\mathbf{s}(t),\mathbf{R}(t)]\in\mathcal{R}_{i}$ and
switches to $u_{a}\left(  t\right)  =j\in\mathcal{N}_{i}$ as soon as the state
vector transitions to a new region $\mathcal{R}_{j}$. Thus, the optimization
problem consists of determining an optimal partition for all $i=1,\ldots,M$
through $g_{i,j}^{\ast}(\mathbf{s}(t),\mathbf{R}(t))$ for all $j\in
\mathcal{N}_{i}$ and the time of \ a transition from $\mathcal{R}_{i}$ to some
$\mathcal{R}_{j}$, $j\in\mathcal{N}_{i}$.

Finally, given control $u_{a}(t)$, the agent's physical dynamics over
$t\in\lbrack t_{a,k},t_{a,k+1})$ are given by
\begin{equation}
\dot{s}_{a}(t)=\left\{
\begin{array}
[c]{cl}%
\frac{X_{i}-s_{a}(t)}{\Vert X_{i}-s_{a}(t)\Vert} & \text{if }t\in\lbrack
t_{a,k}+d_{a,k},t_{a,k+1})\\
0 & \text{otherwise}%
\end{array}
\right.  \label{eq:Dyn_A_param}%
\end{equation}
for some $i=u_{a}(t)\in V$.

\textbf{Parametric control.} As already mentioned, designing an optimal
feedback controller for \textbf{P1} in a 2D space is generally intractable.
The problem remains hard even in the simplified graph topology embedded in the
original 2D space where optimal partitions of the state space must be
determined whenever an agent visits a node. Therefore, an alternative is to
seek a parameterization of these partitions through a parameter vector
$\mathbf{\Theta}$ so as to ultimately replace \textbf{P1} by a problem
requiring the determination\ of an optimal parameter vector $\mathbf{\Theta
}^{\ast}=\arg\min J(\mathbf{\Theta})$ over the set of feasible values of
$\mathbf{\Theta}$. Thus, switching functions of the form $g_{i,j}%
(\mathbf{s}(t),\mathbf{R}(t))$ are expressed as $g_{i,j}(\mathbf{s}%
(t),\mathbf{R}(t);\mathbf{\Theta})$ and an optimal switching function is given
by $g_{i,j}(\mathbf{s}(t),\mathbf{R}(t);\mathbf{\Theta}^{\ast})$.

The parameterization we select in our problem is motivated by the observation
that the movement of agents should be determined based on the values of the
target uncertainty states available to an agent, since the cost function
\eqref{eq:costfunction} is closely related to these values. Thus, we introduce
\emph{threshold parameters} associated with a node $i$ which, when compared to
the actual value of $R_{i}(t)$, provide information about the importance of
visiting this node next when an agent is in its neighborhood and needs to
evaluate the control in (\ref{u_optimal}). We set the thresholds to be
distinct when the agent is at different nodes, thus rendering the control
policy more flexible since it depends on both node uncertainty values and the
agent's position.

\begin{figure}[h]
\centering
\includegraphics[width=0.85\linewidth]{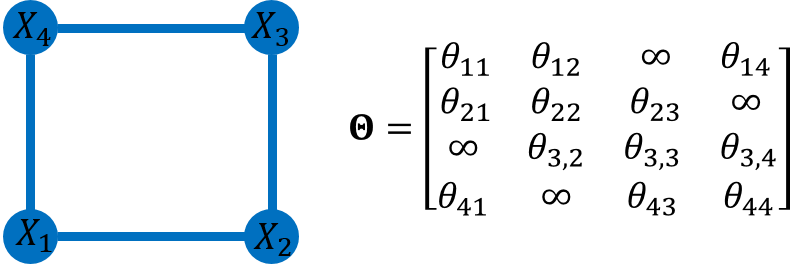}\caption{{\protect\small A
1-agent 4-target example. The target topology graph is shown on the left and
the threshold matrix is on the right.}}%
\label{fig:ProbForm}%
\end{figure}

We represent the node thresholds associated with agent $a$ by an $M\times M$
matrix $\bm\Theta^{a}$ where each row represents the index of the current node
visited by $a$ and a column represents the index of a potential next node to
visit. An example is shown in Fig. \ref{fig:ProbForm} where a threshold
parameter is set to $\infty$ when there is no direct path between the
corresponding nodes. In this example, an agent located at node $1$ uses a
state space partition parametrized by $\theta_{11}$, $\theta_{12}$ and
$\theta_{14}$. The overall parameter matrix accounting for all agents is
denoted by $\bm\Theta$ of dimension $M\times M\times N$.

Next, we define the specific threshold-based controller family we consider.
The starting point is to define a state space region forcing the agent to
remain at node $i$. This is expressed through the condition $R_{i}%
(t)>\theta_{ii}^{a}$. When this is no longer met, i.e., the uncertainty state
at node $i$ is sufficiently low with respect to a level $\theta_{ii}^{a}$,
then the agent may be assigned to a new node $j\neq i$ as long as its
uncertainty state exceeds another threshold, i.e., $R_{j}(t)\geq\theta
_{ij}^{a}$. Since there may be several nodes in the neighborhood of $i$ whose
uncertainty states are high relative to their associated thresholds, we
prioritize nodes in the neighborhood of $i$ by defining an ordered set for
agent $a$ as follows:
\[
\mathcal{N}_{i}^{a}=\{j_{k}\in\mathcal{N}_{i}:j_{1},\ldots,j_{k}%
,\ldots,j_{D_{i}}\}
\]
where $D_{i}$ is the degree of vertex $i$ (the number of edges connected to
vertex $i$). Although the prioritization scheme used may depend on several
factors, in what follows we assume that $\Vert X_{j_{k}}-X_{i}\Vert<\Vert
X_{j_{k+1}}-X_{i}\Vert$ for all $k=1,\ldots,D_{i}$, i.e., the neighbors are
ordered based on their relative proximity to node $i$.

We now define the threshold-based control to specify $u_{a}(t;\bm\Theta)$ in
(\ref{u_optimal}) as follows:
\begin{equation}
u_{a}(t;\bm\Theta)=\left\{
\begin{split}
&  i\qquad\qquad\qquad\qquad\qquad\text{ if }R_{i}(t)>\theta_{ii}^{a}\text{ or
}\\
&  \qquad\qquad\qquad\qquad R_{j}(t)<\theta_{ij}^{a}\text{ for all }%
j\in\mathcal{N}_{i}^{a}\\
\\
&  \arg\min_{\substack{k\\\text{s.t.}j_{k}\in\mathcal{N}_{i}^{a}}}R_{j_{k}%
}\geq\theta_{ij_{k}}^{a}\qquad\text{ otherwise}%
\end{split}
\right.  \label{eq:control_i}%
\end{equation}
Under \eqref{eq:control_i}, the agent first decreases $R_{i}(t)$ below the
threshold $\theta_{ii}^{a}$ before moving to another node in the neighbor set
$\mathcal{N}_{i}^{a}$ with the minimum index $k$ whose associated state
uncertainty value exceeds the threshold $\theta_{ij_{k}}^{a}$. If no such
neighbor exists, the agent remains at the current node maintaining its
uncertainty state under the given threshold level. All agent behaviors are
therefore entirely governed by $\bm\Theta$ through \eqref{eq:control_i}, which
also implicitly determines the dwell time of the agent at node $i$.

\begin{remark}
The controller in \eqref{eq:control_i} is designed to be \emph{distributed} by considering only the states of neighboring nodes and not those of other nodes or of other agents. As such, it is limited to a one-step look-ahead policy. However, it
can be extended to a richer family of more general multi-step look-ahead
policies based on node uncertainty state thresholds. While this causes the
dimensionality of $\bm\Theta$ to increase, the optimization framework
presented in Sec. III is not affected.
\end{remark}

Under \eqref{eq:control_i}, \textbf{P1} is reduced to a simpler parametric
optimization problem of determining the optimal thresholds in matrix
$\bm\Theta^{\star}$ under which the cost function in \eqref{eq:costfunction}
is minimized. Moreover, the resulting agent and node behavior defines a hybrid
system: the node dynamics in (\ref{eq:DynR}) switch between the mode where
$\dot{R}_{i}(t)=0$ and $\dot{R}_{i}(t)=A_{i}-B_{i}N_{i}(t)$ with
$N_{i}(t)=0,1,\ldots,N$, while the agent dynamics in (\ref{eq:Dyn_A_param})
experience a switch whenever there is a sign change in some expression of the
form $(R_{j}(t)-\theta_{ij}^{a})$ as seen in \eqref{eq:control_i}, hence
triggering a control switch. We rewrite the cost in \eqref{eq:costfunction} as
the sum of costs over all intervals $\left[  \tau_{k},\tau_{k+1}\right)  $ for
$k=0,\ldots,K$ where $\tau_{k}$ is the time instant when any of the state
variables experiences a mode switch (these will be explicitly defined as
\textquotedblleft events\textquotedblright\ in the sequel) and $\tau_{K}=T$
denotes the end of the time horizon as shown in \eqref{eq:param_cost}.
Therefore, we have transformed the optimal control problem \textbf{P1} into a
simpler parametric optimization problem \textbf{P2} as follows:
\begin{equation}
\mathbf{P2:}\min_{\bm\Theta\geq\mathbf{0}}\text{ \ }J(\bm\Theta)=\frac{1}%
{T}\sum_{i=1}^{M}\sum_{k=0}^{K}\int_{\tau_{k}}^{\tau_{k+1}}R_{i}(t)\,dt
\label{eq:param_cost}%
\end{equation}
subject to target uncertainty dynamics \eqref{eq:DynR}, agent state dynamics
(\ref{eq:Dyn_A_param}) and the control policy \eqref{eq:control_i}.

\begin{remark}
Using the optimal threshold matrix $\bm\Theta^{\star}$, the
optimal dwell times and target visiting sequences can both be determined on
line while executing the control policy \eqref{eq:control_i}. It is
interesting to note that, despite the a priori prioritization imposed in
$\mathcal{N}_{i}^{a}$, the actual target visiting sequence will be adjusted as
a result of the thresholds being adjusted during the optimization process.
This is because the optimization process will decrease the threshold values of
nodes that maintain higher uncertainties, hence inducing agents to visit them more frequently.
\end{remark}

\section{Infinitesimal Perturbation Analysis (IPA)}

\label{sec:IPA} In the previous section, agent trajectories are selected from
the family $\mathbf{s}(\bm\Theta,\mathbf{s}_{0},\mathbf{R}_{0})$ with
parameter $\bm\Theta$ and given initial agent positions $\mathbf{s}_{0}$ and
node uncertainty states $\mathbf{R}_{0}$. The state dynamics are governed by
\eqref{eq:DynR} and (\ref{eq:Dyn_A_param}) under the control policy
\eqref{eq:control_i}. An \textquotedblleft event\textquotedblright\ is defined
as any discontinuous change in any one of the state variables (e.g., a
threshold has been met by some $R_{i}(t)$). The $k$-th event occurrence time
is denoted by $\tau_{k}(\bm\Theta)$. We use Infinitesimal Perturbation
Analysis (IPA) to obtain on line the gradient of the cost function in
(\ref{eq:param_cost}) with respect to the parameters in $\bm\Theta$, hence
seeking an optimal solution through a gradient descent process. IPA specifies
how changes in the parameter $\bm\Theta$ influence event times $\tau
_{k}(\bm\Theta)$, $k=1,2,\ldots$, the trajectories $\mathbf{s}(\bm\Theta
,\mathbf{s}_{0},\mathbf{R}_{0})$, and ultimately the cost function
\eqref{eq:param_cost}. We first briefly review the IPA framework for general
hybrid systems as presented in \cite{cassandras2010perturbation} and then
apply it to our specific setting.

Let $\{\tau_{k}(\theta)\}$, $k=1,\ldots,K$, denote the occurrence times of all
events in the state trajectory of a hybrid system with dynamics $\dot
{x}\ =\ f_{k}(x,\theta,t)$ over an interval $[\tau_{k}(\theta),\tau
_{k+1}(\theta))$, where $\theta\in\Theta$ is some parameter vector and
$\Theta$ is a given compact, convex set. For convenience, we set $\tau_{0}=0$
and $\tau_{K+1}=T$. We use the Jacobian matrix notation: $x^{\prime}%
(t)\equiv\frac{\partial x(\theta,t)}{\partial\theta}$ and $\tau_{k}^{\prime
}\equiv\frac{\partial\tau_{k}(\theta)}{\partial\theta}$, for all state and
event time derivatives. It is shown in \cite{cassandras2010perturbation} that
\begin{equation}
\frac{d}{dt}x^{\prime}(t)=\frac{\partial f_{k}(t)}{\partial x}x^{\prime
}(t)+\frac{\partial f_{k}(t)}{\partial\theta}, \label{eq:IPA_1}%
\end{equation}
for $t\in\lbrack\tau_{k},\tau_{k+1})$ with boundary condition:
\begin{equation}
x^{\prime}(\tau_{k}^{+})=x^{\prime}(\tau_{k}^{-})+[f_{k-1}(\tau_{k}^{-}%
)-f_{k}(\tau_{k}^{+})]\tau_{k}^{\prime} \label{eq:IPA_2}%
\end{equation}
for $k=1,...,K$. In order to complete the evaluation of $x^{\prime}(\tau
_{k}^{+})$ in (\ref{eq:IPA_2}), we need to determine $\tau_{k}^{\prime}$. If
the event at $\tau_{k}$ is \emph{exogenous} (i.e., independent of $\theta$),
$\tau_{k}^{\prime}=0$. However, if the event is \emph{endogenous}, there
exists a continuously differentiable guard function $g_{k}:\mathbb{R}%
^{n}\times\Theta\rightarrow\mathbb{R}$ such that $\tau_{k}\ =\ \min
\{t>\tau_{k-1}\ :\ g_{k}\left(  x\left(  \theta,t\right)  ,\theta\right)
=0\}$ and
\begin{equation}
\tau_{k}^{\prime}=-[\frac{\partial g_{k}}{\partial x}f_{k}(\tau_{k}^{-}%
)]^{-1}(\frac{\partial g_{k}}{\partial\theta}+\frac{\partial g_{k}}{\partial
x}x^{\prime}(\tau_{k}^{-})) \label{eq:IPA_3}%
\end{equation}
as long as $\frac{\partial g_{k}}{\partial x}f_{k}(\tau_{k}^{-})\neq0$
(details can be found in \cite{cassandras2010perturbation}).

Differentiating the cost $J(\bm\Theta)$ in \textbf{P2}, we obtain
\begin{align}
\nabla J(\bm\Theta)=  &  \frac{1}{T}\sum_{i=1}^{M}\sum_{k=0}^{K}%
\Bigg(\int_{\tau_{k}}^{\tau_{k+1}}\nabla R_{i}(t)\,dt\nonumber\\
&  +R_{i}(\tau_{k+1})\nabla\tau_{k+1}-R_{i}(\tau_{k})\nabla\tau_{k}%
\Bigg)\nonumber\\
=  &  \frac{1}{T}\sum_{i=1}^{M}\sum_{k=0}^{K}\int_{\tau_{k}}^{\tau_{k+1}%
}\nabla R_{i}(t)\,dt \label{eq:gradient_J}%
\end{align}
where the gradient operator $\nabla\equiv\frac{\partial}{\partial\bm\Theta}$
and all terms of the form $R_{i}(\tau_{k})\nabla\tau_{k}$ for all $k$ are
cancelled with $\tau_{0}=0$ and $\tau_{K}=T$ fixed. We first derive the
integrand $\nabla R_{i}(t)$ in \eqref{eq:gradient_J} for $i=1,\ldots,M$ and
then integrate over $[0,T]$ to obtain $\nabla J(\bm\Theta)$. The following
lemma shows that the integrand $\nabla R_{i}(t)$ remains constant between
any two consecutive events and can be updated only at some event time. This
establishes the fully \emph{event-driven} nature of our IPA-based gradient algorithm.

\begin{lemma}
$\nabla R_{i}(t)$ remains constant for $t \in[\tau_{k}, \tau_{k+1} ), k =
0,\ldots K-1$. \label{lemma:lemma1}
\end{lemma}

\textbf{Proof.}
In each inter-event interval, $\dot{R}_{i}(t)$ in \eqref{eq:DynR} remains
constant and, therefore, $\frac{\partial f_{k}(t)}{\partial R_{i}} = 0$ and
$\frac{\partial f_{k}(t)}{\partial\theta} = 0$ where either $f_{k}%
(t)=A_{i}-B_{i}N_{i}(t)$ or $f_{k}(t)=0$. From \eqref{eq:IPA_1}, we can obtain
$\frac{d}{dt} R^{\prime}_{i} = 0$. As a result,
\begin{equation}
\nabla R_{i}(t) = \nabla R_{i}(\tau_{k}^{+}), \quad t \in[\tau_{k}, \tau_{k+1}
)
\end{equation}
$\blacksquare$

\begin{figure}[h]
\centering
\includegraphics[width=1.0\linewidth]{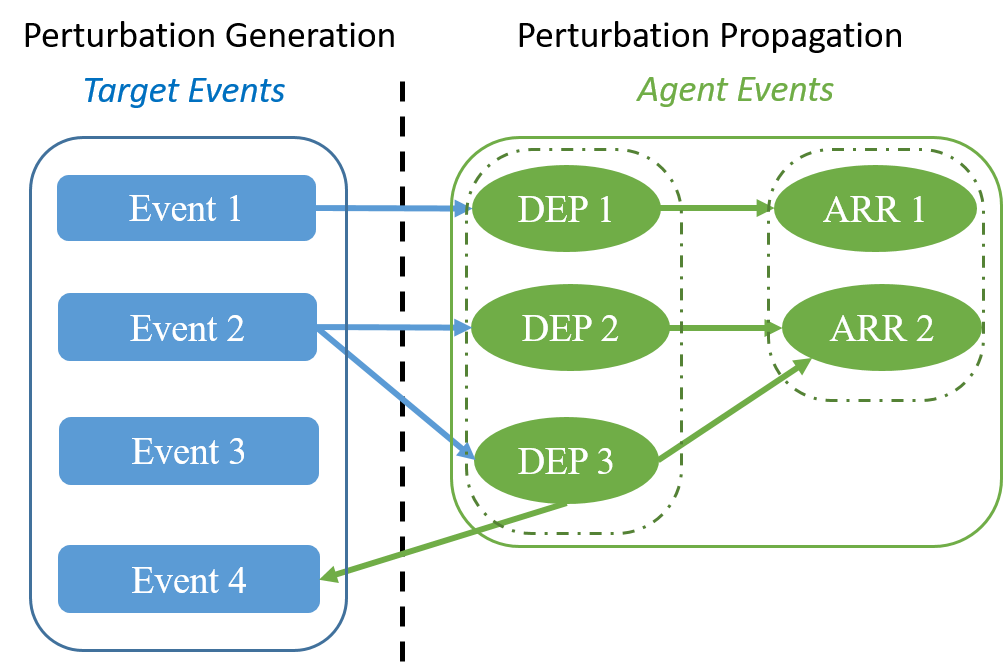}\caption{{\protect\small Event
inducing scheme and the corresponding process of perturbation generation and
propagation.}}%
\label{fig:pert_prop}%
\end{figure}

In the following, we will show the derivation of $\nabla R_{i}(t)$ at each
event time $\tau_{k}$. To do so, we need to first define all events in this
hybrid system which may cause discontinuities in $\nabla R_{i}(t)$. In view of
{\eqref{eq:control_i}}, there are four types of \textquotedblleft target
events\textquotedblright\ (labeled Event $1$ to $4$) corresponding to
$R_{i}(t)$ crossing some threshold value from above/below, reaching the value
$R_{i}(t)=0$ from above or leaving the value $R_{i}(t)=0$. In our parametric
control problem \textbf{P2}, each agent's movement is controlled by the target
thresholds. Therefore, a target event may induce a switch in the agent
dynamics \eqref{eq:DynR} through an agent departure event, denoted by
\textbf{DEP}, occurring at $t_{a,k}+d_{a,k}$ in (\ref{eq:Dyn_A_param}). In
turn, this event will induce this agent's arrival event at the next node
visited, denoted by \textbf{ARR}. The process of how events can induce other
events is detailed next and is graphically summarized in Fig.
\ref{fig:pert_prop}.

For notational simplicity, we use $\downarrow=$ as an operator indicating that
the value on its left-hand-side reaches the value on its right-hand-side from
above. Similarly, $\uparrow=$ means reaching from below, and $=\uparrow$ means
increasing from the value on the right-hand-side. In addition,
since the derivative with respect to $\bm\Theta$ is updated differently at
different entries, we use $p,q,z$ to indicate the $pq$-entry of the parameter
$\bm\Theta$ of agent $z$. 

\textbf{Event 1: $R_{i}(\tau_{k})\downarrow=\theta_{ii}^{a}$ }. In this case,
$R_{i}(t)$ reaches the threshold $\theta_{ii}^{a}$ from above. It is an
endogenous event and the guard condition is $g_{k}=R_{i}-\theta_{ii}^{a}=0$ in
\eqref{eq:IPA_3}. Therefore, the event time derivative with respect to the
$pq$-th entry of the parameter $\bm\Theta$ of agent $z$ is as follows:

\begin{equation}
(\tau_{k}^{\prime})_{pq}^{z}=\left\{
\begin{array}
[c]{ll}%
-\frac{-1+\left(  R_{i}^{\prime}(\tau_{k}^{-})\right)  _{pq}^{z}}{A_{i}%
-B_{i}N_{i}(\tau_{k}^{-})} & \text{if }p=q=i,\text{ and }z=a\\
-\frac{\left(  R_{i}^{\prime}(\tau_{k}^{-})\right)  _{pq}^{z}}{A_{i}%
-B_{i}N_{i}(\tau_{k}^{-})} & \text{otherwise }%
\end{array}
\right.  \label{eq:event1}%
\end{equation}
Based on \eqref{eq:control_i}, this event may induce an agent departure from
its current node location which we denote as event DEP1. Through this event,
the value of the event time derivative in \eqref{eq:event1} will be
transferred to $R_{i}^{\prime}(t)$ as shown next.

\textbf{DEP1: Agent departure event 1.} In this case, the agent departure is
induced by Event $1$. Using \eqref{eq:IPA_2} and \eqref{eq:event1}, we obtain
\begin{align}
&  \left(  R_{i}^{\prime}(\tau_{k}^{+})\right)  _{pq}^{z}=\left(
R_{i}^{\prime}(\tau_{k}^{-})\right)  _{pq}^{z}-B_{i}\left(  \tau_{k}^{\prime
}\right)  _{pq}^{z}\nonumber\\
&  =\left\{
\begin{array}
[c]{ll}%
\frac{A_{i}-B_{i}\left(  N_{i}(\tau_{k}^{-})-1\right)  }{A_{i}-B_{i}N_{i}%
(\tau_{k}^{-})}\left(  R_{i}^{\prime}(\tau_{k}^{-})\right)  _{pq}^{z}%
-\frac{B_{i}}{A_{i}-B_{i}N_{i}(\tau_{k}^{-})} & \\
\qquad\qquad\qquad\qquad\qquad\quad\text{if }p=q=i,\text{ and }z=a & \\
\frac{A_{i}-B_{i}\left(  N_{i}(\tau_{k}^{-})-1\right)  }{A_{i}-B_{i}N_{i}%
(\tau_{k}^{-})}\left(  R_{i}^{\prime}(\tau_{k}^{-})\right)  _{pq}^{z}%
\quad\text{otherwise} &
\end{array}
\right.  \label{eq:event5_1}%
\end{align}
This agent departure event will eventually induce an arrival event at another
target. The value of the event time derivative in \eqref{eq:event1} will be
transferred to this arrival event, therefore, $\tau_{k+1}^{\prime}=\tau
_{k}^{\prime}$ because the travel time between any two nodes $i$ and $j$ is
fixed and independent of $\bm\Theta$. To see this, set $g_{k+1}=g_{k}+c$ where
$c$ is a constant determined by the travel time. Through \eqref{eq:IPA_3}, it
is obvious that adding a constant after $g_{k}$ does not affect the
derivative. Therefore, we can transfer the value of $\tau_{k}^{\prime}$ to
$\tau_{k+1}^{\prime}$ (a similar proof can be found in
\cite{cassandras2010perturbation} Lemma 2.1).

\textbf{ARR1: Agent arrival event 1.} This is induced by the earlier DEP1 at
node $i$, which is again induced by the target event $R_{i}(\tau
_{k})\downarrow=\theta_{ii}^{a}$ (Event 1) and we transfer the value of the
event time derivative to obtain
\begin{equation}
\left(  \tau_{k+1}^{\prime}\right) _{pq}^{z}=\left(  \tau_{k}^{\prime}\right)
_{pq}^{z}=\left\{
\begin{array}
[c]{ll}%
-\frac{-1+\left(  R_{i}^{\prime}(\tau_{k}^{-})\right) _{pq}^{z}}{A_{i}%
-B_{i}N_{i}(\tau_{k}^{-})}\quad\text{if }p=q=i, & \\
\qquad\qquad\qquad\qquad\text{ and }z=a & \\
-\frac{\left( R_{i}^{\prime}(\tau_{k}^{-})\right) _{pq}^{z}}{A_{i}-B_{i}%
N_{i}(\tau_{k}^{-})}\quad\text{otherwise} &
\end{array}
\right.
\end{equation}
and through \eqref{eq:IPA_2},
\begin{equation}%
\begin{split}
&  \left(  R_{j}^{\prime}(\tau_{k+1}^{+})\right)  _{pq}^{z}=\left(
R_{j}^{\prime}(\tau_{k+1}^{-})\right)  _{pq}^{z}+B_{j}\left(  \tau_{k}%
^{\prime}\right)  _{pq}^{z}\\
&  =\hspace{-1mm}\left\{
\begin{array}
[c]{ll}%
\left(  R_{j}^{\prime}(\tau_{k+1}^{-})\right)  _{pq}^{z}-\frac{B_{j}}%
{A_{i}-B_{i}N_{i}(\tau_{k}^{-})}\left(  \left(  R_{i}^{\prime}(\tau_{k}%
^{-})\right)  _{pq}^{z}-1\right)  & \\
\qquad\qquad\qquad\qquad\qquad\qquad\text{ if }p=q=i,\text{ and }z=a & \\
\left(  R_{j}^{\prime}(\tau_{k+1}^{-})\right) _{pq}^{z}-\frac{B_{j}}%
{A_{i}-B_{i}N_{i}(\tau_{k}^{-})}\left(  R_{i}^{\prime}(\tau_{k}^{-})\right)
_{pq}^{z}\quad\text{otherwise} &
\end{array}
\right.
\end{split}
\label{eq:event6_1}%
\end{equation}

\textbf{Event 2: $R_{j}(\tau_{k})\uparrow=\theta_{ij}^{a}$ }. This event
occurs when an agent is at node $i$ and $R_{j}(t)$ at $j\neq i$ exceeds the
threshold $\theta_{ij}^{a}$. The event is endogenous and the guard condition
in \eqref{eq:IPA_3} is $g_{k}=R_{j}-\theta_{ij}^{a}=0$. The event time
derivative is obtained from \eqref{eq:IPA_3} as follows:
\begin{equation}
\left(  \tau_{k}^{\prime}\right)  _{pq}^{z}=\left\{
\begin{array}
[c]{ll}%
-\frac{-1+\left(  R_{j}^{\prime}(\tau_{k}^{-})\right)  _{pq}^{z}}{A_{j}%
-B_{j}N_{j}(\tau_{k}^{-})} & \text{if }p=i,q=j,\text{and }z=a\\
-\frac{\left(  R_{j}^{\prime}(\tau_{k}^{-})\right)  _{pq}^{z}}{A_{j}%
-B_{j}N_{j}(\tau_{k}^{-})} & \text{otherwise }%
\end{array}
\right.  \label{eq:event2}%
\end{equation}
Looking at \eqref{eq:control_i}, this event can induce an agent departure
event depending on whether $R_{i}(\tau_{k})>0$ or not: in the former case, the
event is denoted by DEP2 and in the latter it is denoted by DEP3. The value of
the derivative in \eqref{eq:event2} will be transferred to $R_{i}^{\prime}(t)$
through one of these agent departure events.

\textbf{DEP2: Agent departure event 2.} In this case, $R_{i}(\tau_{k})>0$.
Using \eqref{eq:IPA_2} and the event time derivative in \eqref{eq:event2}, we
obtain
\begin{align}
&  \left(  R_{i}^{\prime}(\tau_{k}^{+})\right)  _{pq}^{z}=\left(
R_{i}^{\prime}(\tau_{k}^{-})\right)  _{pq}^{z}-B_{i}\left(  \tau_{k}^{\prime
}\right)  _{pq}^{z}\nonumber\\
&  =\left\{
\begin{array}
[c]{ll}%
\left(  R_{i}^{\prime}(\tau_{k}^{-})\right)  _{pq}^{z}+\frac{B_{i}}%
{A_{j}-B_{j}N_{j}(\tau_{k}^{-})}\left(  \left(  R_{i}^{\prime}(\tau_{k}%
^{-})\right)  _{pq}^{z}-1\right)  & \\
\qquad\qquad\qquad\qquad\qquad\quad\text{ if }p=i, q=j,\text{and }z=a & \\
\left(  R_{i}^{\prime}(\tau_{k}^{-})\right)  _{pq}^{z}+\frac{B_{i}}%
{A_{j}-B_{j}N_{j}(\tau_{k}^{-})}\left(  R_{i}^{\prime}(\tau_{k}^{-})\right)
_{pq}^{z}\quad\text{otherwise } &
\end{array}
\right.  \label{eq:event5_2}%
\end{align}

\textbf{DEP3: Agent departure event 3.} This event is complementary to DEP2
where the agent departure is induced by Event $2$ but $R_{i}(\tau_{k})=0$.
Based on (\ref{eq:DynR}), the target dynamics after this event either remain
$\dot{R}_{i}(t)=0$ or switch to $\dot{R}_{i}(t)=A_{i}-B_{i}N_{i}(\tau_{k}%
^{+})$ depending on whether $A_{i}>B_{i}N_{i}(\tau_{k}^{+})$ or not. Thus,
there are two sub-cases to consider as follows.

\textbf{DEP3-1:} $A_{i}>B_{i}N_{i}(\tau_{k}^{+})$. In this sub-case, the
target dynamics switch from $\dot{R}_{i}(t)=0$ for $t\in\lbrack\tau_{k-1}%
,\tau_{k})$ to $\dot{R}_{i}(t)=A_{i}-B_{i}N_{i}(t)$ for $t\in\lbrack\tau
_{k},\tau_{k+1})$. We know $R_{i}^{\prime}(\tau_{k}^{-})=0$ because
$R_{i}(\tau_{k})=0$ before the agent departure and the value $R_{i}^{\prime
}(t)=0$ holds as long as $R_{i}(t)=0$. Using \eqref{eq:IPA_2} and the event
time derivative in \eqref{eq:event2}, we obtain
\begin{align}
&  \left(  R_{i}^{\prime}(\tau_{k}^{+})\right)  _{pq}^{z}=-(A_{i}-B_{i}%
N_{i}(\tau_{k}^{+}))\left(  \tau_{k}^{\prime}\right)  _{pq}^{z}\nonumber\\
&  =\left\{
\begin{array}
[c]{ll}%
\frac{A_{i}-B_{i}N_{i}(\tau_{k}^{+})}{A_{j}-B_{j}N_{j}(\tau_{k}^{-})}\left(
\left(  R_{j}^{\prime}(\tau_{k}^{-})\right)  _{pq}^{z}-1\right)  \quad\text{if
}p=i,q=j, & \\
\qquad\qquad\qquad\qquad\qquad\qquad\qquad\quad\text{ and }z=a & \\
\frac{A_{i}-B_{i}N_{i}(\tau_{k}^{+})}{A_{j}-B_{j}N_{j}(\tau_{k}^{-})}\left(
R_{j}^{\prime}(\tau_{k}^{-})\right)  _{pq}^{z}\qquad\qquad\text{otherwise } &
\end{array}
\right.  \label{eq:event5_3}%
\end{align}

\textbf{DEP3-2:} $A_{i}\leq B_{i}N_{i}(\tau_{k}^{+})$. In this sub-case, the
target dynamics remain $\dot{R}_{i}(t)=0$ before and after the event at
$\tau_{k}$. Therefore, the state dynamics in \eqref{eq:IPA_2} satisfy
$f_{k-1}(\tau_{k}^{-})=f_{k}(\tau_{k}^{+})$ and we have
\begin{equation}
R_{i}^{\prime}(\tau_{k}^{+})=0\quad\text{for all } p,q,z\\
\end{equation}

\begin{remark}
Note that DEP3-1 induces another target event (Event $4$)
since $R_{i}(t)$ increases after the agent's departure. Moreover, both DEP2
and DEP3 will induce an agent arrival event at the next visiting target.
\end{remark}

\textbf{ARR2: Agent arrival event 2.} This event is induced by an earlier
agent departure event at a target $i$ which is again induced by the previous
Event 2 $R_{j}(\tau_{k})\uparrow=\theta_{ij}^{a}$. Similar to the derivation
in \textbf{ARR1}, we transfer the prior event time derivative value to the
current arrival time derivative:
\begin{equation}
\left(  \tau_{k+1}^{\prime}\right)  _{pq}^{z}=\left(  \tau_{k}^{\prime
}\right)  _{pq}^{z}=\left\{
\begin{array}
[c]{ll}%
-\frac{-1+\left(  R_{j}^{\prime}(\tau_{k}^{-})\right)  _{pq}^{z}}{A_{j}%
-B_{j}N_{j}(\tau_{k}^{-})}\quad\text{if }p=i,q=j, & \\
\qquad\qquad\qquad\qquad\quad\text{ and }z=a & \\
-\frac{\left(  R_{j}^{\prime}(\tau_{k}^{-})\right)  _{pq}^{z}}{A_{j}%
-B_{j}N_{j}(\tau_{k}^{-})}\qquad\text{otherwise } &
\end{array}
\right.
\end{equation}
Through \eqref{eq:IPA_2} we obtain
\begin{equation}%
\begin{split}
&  \left(  R_{j}^{\prime}(\tau_{k+1}^{+})\right)  _{pq}^{z}=\left(
R_{j}^{\prime}(\tau_{k+1}^{-})\right)  _{pq}^{z}+B_{j}\left(  \tau_{k}%
^{\prime}\right)  _{pq}^{z}\\
&  =\left\{
\begin{array}
[c]{ll}%
\left(  R_{j}^{\prime}(\tau_{k+1}^{-})\right)  _{pq}^{z}-\frac{B_{j}}%
{A_{j}-B_{j}N_{j}(\tau_{k}^{-})}\left(  \left(  R_{j}^{\prime}(\tau_{k}%
^{-})\right)  _{pq}^{z}-1\right)  & \\
\qquad\qquad\qquad\qquad\qquad\quad\text{ if }p=i,q=j,\text{ and }z=a & \\
\left(  R_{j}^{\prime}(\tau_{k+1}^{-})\right)  _{pq}^{z}-\frac{B_{j}}%
{A_{j}-B_{j}N_{j}(\tau_{k}^{-})}\left(  R_{j}^{\prime}(\tau_{k}^{-})\right)
_{pq}^{z}\text{ otherwise} &
\end{array}
\right.
\end{split}
\label{eq:event6_2}%
\end{equation}
Notice here that $\tau_{k}$ is the prior agent departure time and $\tau_{k+1}$
is the current agent arrival time and the derivatives $R_{j}^{\prime}%
(\tau_{k+1}^{-})$ and $R_{j}^{\prime}(\tau_{k}^{-})$ can be different since
$R_{j}^{\prime}(\tau_{k}^{-})$ may change due to arrivals or departures of
other agents during $[\tau_{k},\tau_{k+1})$.

\textbf{Event 3:} $R_{i}(t)\downarrow=0$. This event corresponds to the target
uncertainty state reaching zero from above, therefore from (\ref{eq:DynR}) the
target state dynamics switch from $\dot{R}_{i}(t)=A_{i}-B_{i}N_{i}(t),$
$t\in\lbrack\tau_{k-1},\tau_{k})$ to $\dot{R}_{i}(t)=0,$ $t\in\lbrack\tau
_{k},\tau_{k+1})$. It is an endogenous event that occurs when $g_{k}%
(x,\theta)=R_{i}=0$. According to \eqref{eq:IPA_3},
\begin{equation}
\left(  \tau_{k}^{\prime}\right)  _{pq}^{z}=-\frac{\left(  R_{i}^{\prime}%
(\tau_{k}^{-})\right)  _{pq}^{z}}{A_{i}-B_{i}N_{i}(\tau_{k}^{-})}
\label{eq:event3_1}%
\end{equation}
Replacing $\tau_{k}^{\prime}$ in \eqref{eq:IPA_2} with the result in
\eqref{eq:event3_1}, we have
\begin{align}
\left(  R_{i}^{\prime}(\tau_{k}^{+})\right)  _{pq}^{z}=  &  \left(
R_{i}^{\prime}(\tau_{k}^{-})\right)  _{pq}^{z}+\left(  A_{i}-B_{i}N_{i}%
(\tau_{k}^{-})\right)  \left(  \tau_{k}^{\prime}\right)  _{pq}^{z}%
=0\nonumber\\
&  \qquad\qquad\qquad\qquad\qquad\text{ for all }p,q,z
\end{align}
This indicates that $\nabla R_{i}(t)$ is always reset to $0$ whenever the
target's uncertainty state is reduced to zero. This is an uncontrollable event
and does not induce any other event.

\textbf{Event 4:} $R_{i}(t)=\uparrow0$. In this case, the target value leaves
zero and the dynamics in (\ref{eq:DynR}) switch from $\dot{R}_{i}(t)=0,$
$t<\tau_{k}$ to $\dot{R}_{i}(t)=A_{i}-B_{i}N_{i}(t),$ $t\geq\tau_{k}$. This
event is induced by an agent departure event (DEP3) which is in turn induced
by Event $2$. This is an exogenous event functioning only as an indicator of
$R_{i}(t)$ increasing from zero. Therefore, $\tau_{k}^{\prime}= 0$ and the
derivative $R_{i}^{\prime}(t)$ will not be affected.

\begin{remark}
The analysis of Events $1$ to $4$ shows that all non-zero
gradient values are caused by target events and then propagated through the
various agent departure and arrival events.
\end{remark}

\textbf{IPA-based gradient descent algorithm} Once we have derived the
gradient $\nabla J(\bm\Theta)$ in \eqref{eq:gradient_J}, we update the
parameter $\bm\Theta$ based on a standard gradient descent scheme as follows.
\begin{equation}
\bm\Theta^{(l+1)}=\Pi\left[  \bm\Theta^{(l)}-\beta^{l}\nabla J(\bm\Theta
^{(l)})\right]  \label{eq:param_update}%
\end{equation}
where the operator $\Pi\equiv\max\{\bm\cdot,\bm0\}$, $l$ indexes the number of
iterations, and $\beta^{l}$ is a diminishing step-size sequence satisfying
$\sum_{l=0}^{\infty}\beta^{l}=\infty$ and $\lim_{l\rightarrow\infty}\beta
^{l}=0$.

\section{One-agent case analysis}

\label{sec:One-agent case analysis} Recalling our control policy in
\eqref{eq:control_i}, the diagonal entries in the parameter matrix control the
dwell times at nodes, whereas the off-diagonal entries control the feasible
node visiting sequence. In what follows, we will show that in a single-agent
case the optimal values of diagonal entries in
\eqref{eq:single_agent_param_table} are always zero. This structural property
indicates that the agent visiting a node should always reduce the uncertainty
state to zero before moving to the next node.

Ignoring the superscript agent index, the single-agent threshold matrix is
written as
\begin{equation}
\bm\Theta=%
\begin{bmatrix}
\theta_{11} & \theta_{12} & \theta_{13} & \dots & \theta_{1M}\\
\theta_{21} & \theta_{22} & \theta_{23} & \dots & \theta_{2M}\\
\vdots & \vdots & \vdots & \ddots & \vdots\\
\theta_{M1} & \theta_{M2} & \theta_{M3} & \dots & \theta_{MM}%
\end{bmatrix}
\label{eq:single_agent_param_table}%
\end{equation}

\begin{assumption}
For any $\epsilon>0$, there exists a finite time horizon $T>t_{K}-\frac
{c}{1-\epsilon}$ where $t_{K}$ is such that $\Vert\nabla R_{i}(t_{1})-\nabla R_{i}(t_{2})\Vert\leq\epsilon/M$, $i = 1,\ldots,M$ for all $t_{1},t_{2}>t_{K}$ and $c$ is a
finite constant.
\end{assumption}

\begin{assumption}
The current node visiting sequence is optimal.
\end{assumption}

The first assumption is a technical one and it ensures that the optimization
problem is defined over a sufficiently long time horizon $T$ to allow the
gradient to converge. The second assumption allows us to reduce the parameter
matrix \eqref{eq:single_agent_param_table} to a vector of its diagonal
elements only:
\begin{equation}
\Theta_{d}=\left[  \theta_{1},\theta_{2},\ldots,\theta_{M}\right]  ^{\top}%
\geq\mathbf{0}_{M\times1}. \label{eq:1A_param_simp}%
\end{equation}

\begin{theorem}
Consider $M$ targets and a single agent under the parametric control
$\Theta_{d}$. The optimal thresholds satisfy $\Theta_{d}^{\star}%
=\mathbf{0}_{M\times1}$. \label{thm:thm1}
\end{theorem}

\textbf{Proof.}
To establish the proof, we will show that the derivative $\partial
J(\Theta_{d})/\partial\theta_{i}$ satisfies $\partial J(\Theta_{d}%
)/\partial\theta_{i}>0$ for every $i=1,\ldots,M$. As a result, through the
parameter update scheme \eqref{eq:param_update}, $\Theta_{d}$ will eventually
be reduced to $\mathbf{0}$. First, for every element $\theta_{i}$ in
$\Theta_{d}$, we have
\begin{equation}
\frac{\partial J(\Theta_{d})}{\partial\theta_{i}}=\frac{1}{T}\int_{t=0}%
^{T}\sum_{m=1}^{M}\frac{\partial R_{m}}{\partial\theta_{i}}%
\,dt\label{eq:thm1_1}%
\end{equation}
The value of the integrand over time is given by the IPA results in Sec. III
as follows:

\textbf{Agent departures.} In the single-agent case, all agent departure
events are of type DEP1. From \eqref{eq:event5_1}, the IPA derivatives with
respect to each element in $\Theta_{d}$ after such events are:
\begin{equation}
\left\{
\begin{split}
&  \frac{\partial R_{i}}{\partial\theta_{i}}(\tau_{k}^{+})=\frac{A_{i}}%
{A_{i}-B_{i}}\frac{\partial R_{i}}{\partial\theta_{i}}(\tau_{k}^{-}%
)-\frac{B_{i}}{A_{i}-B_{i}}\\
&  \frac{\partial R_{i}}{\partial\theta_{j}}(\tau_{k}^{+})=\frac{A_{i}}%
{A_{i}-B_{i}}\frac{\partial R_{i}}{\partial\theta_{j}}(\tau_{k}^{-}%
)\quad\text{ for }j\neq i
\end{split}
\right.  \label{eq:thm1_3}%
\end{equation}
\textbf{Agent arrivals.} An agent arrival event at node $i$ is induced by the
earlier DEP1 event at some previously visited node $j\neq i$. According to
\eqref{eq:event6_1}, the IPA derivatives are:
\begin{equation}
\left\{
\begin{split}
&  \frac{\partial R_{i}}{\partial\theta_{j}}(\tau_{k}^{+})=\frac{\partial
R_{i}}{\partial\theta_{j}}(\tau_{k}^{-})-\frac{B_{i}}{A_{j}-B_{j}}\left(
\frac{\partial R_{j}}{\partial\theta_{j}}(\tau_{k}^{-})-1\right)  \\
&  \qquad\qquad\qquad\qquad\qquad\qquad\text{where $j$ is the prior node}\\
&  \frac{\partial R_{i}}{\partial\theta_{l}}(\tau_{k}^{+})=\frac{\partial
R_{i}}{\partial\theta_{l}}(\tau_{k}^{-})-\frac{B_{i}}{A_{j}-B_{j}}%
\frac{\partial R_{j}}{\partial\theta_{l}}(\tau_{k}^{-})\text{ for }l\neq j
\end{split}
\right.  \label{eq:thm1_4}%
\end{equation}
To simplify the notation, we set
\begin{equation}
\nabla R_{i}(t)=\left[  \frac{\partial R_{i}}{\partial\theta_{1}}%
,\frac{\partial R_{i}}{\partial\theta_{2}},\ldots,\frac{\partial R_{i}%
}{\partial\theta_{M}}\right]  ^{\top}%
\end{equation}
and
\begin{equation}
\nabla R(t)=\left[  \nabla R_{1},\nabla R_{2},\ldots,\nabla R_{M}\right]
^{\top}.\label{eq:R_param_simp}%
\end{equation}
The evolution of the gradient vector $\nabla R(t)$ follows a system of linear
equations in \eqref{eq:thm1_3} and \eqref{eq:thm1_4} for each node $i$.
Solving this system of equations, we obtain the only possible equilibrium for
every node $i=1,\dots,M$:
\begin{equation}
\frac{\partial R_{i}}{\partial\theta_{i}}=1\text{ and }\frac{\partial R_{i}%
}{\partial\theta_{j}}=0\text{ for }j\neq i\label{eq:thm1_eq}%
\end{equation}
Using Assumption 1, for any $0<\epsilon<1$, there exists a $t_{K}$ such that
$\Vert\frac{\partial R_{i}}{\partial\theta_{i}}(t_{K})-1\Vert<\epsilon/M$ for
all $i=1,\ldots,M$ and $\Vert\frac{\partial R_{i}}{\partial\theta_{j}}%
(t_{K})\Vert<\epsilon/M$ for all $i\neq j$. We now rewrite \eqref{eq:thm1_1}
with the integral separated into two parts over $[0,T]$ as follows:
\[
\frac{\partial J(\Theta_{d})}{\partial\theta_{i}}=\frac{1}{T}\left(
\int_{t=0}^{t_{K}}\sum_{m=1}^{M}\frac{\partial R_{m}(t)}{\partial\theta_{i}%
}dt+\int_{t_{K}}^{T}\sum_{m=1}^{M}\frac{\partial R_{m}(t)}{\partial\theta_{i}%
}dt\right)
\]
The first integral above corresponds to the transient stage before $t_{K}$ and
there exists some constant $c$ whose value is smaller that this integral so
that
\begin{equation}%
\begin{split}
\frac{\partial J(\Theta_{d})}{\partial\theta_{i}} &  \geq\frac{1}{T}\left(
c+\int_{t_{K}}^{T}\sum_{m=1}^{M}\frac{\partial R_{m}(t)}{\partial\theta_{i}%
}dt\right)  \\
&  =\frac{1}{T}\left(  c+\int_{t_{K}}^{T}\left[  \frac{\partial R_{i}%
(t)}{\partial\theta_{i}}+\sum_{m\neq i}\frac{\partial R_{m}(t)}{\partial
\theta_{i}}\right]  dt\right)  \\
& \hspace{-15mm} =\frac{1}{T}\left(  c+\int_{t_{K}}^{T}dt+\int_{t_{K}}^{T}\left[
\frac{\partial R_{i}(t)}{\partial\theta_{i}}-1+\sum_{m\neq i}\frac{\partial
R_{m}(t)}{\partial\theta_{i}}\right]  dt\right)  \\
&  \hspace{-15mm}\geq\frac{1}{T}\left(  c\!+\!T\!-\!t_{K}-\int_{t_{K}}%
^{T}\!\Vert\frac{\partial R_{i}(t)}{\partial\theta_{i}}-1\Vert+\sum_{m\neq
i}\Vert\frac{\partial R_{m}(t)}{\partial\theta_{i}}\Vert dt\right)  \\
&  \hspace{-15mm}\geq\frac{1}{T}\left(  c+T-t_{K}-\int_{t_{K}}^{T}%
\epsilon/M+(M-1)\epsilon/M\,dt\right)  \\
&  \hspace{-15mm}=\frac{1}{T}\left(  c+\left(  1-\epsilon\right)
(T-t_{K})\right)
\end{split}
\label{eq:thm1_5}%
\end{equation}
Therefore, as long as $T>t_{K}-\frac{c}{1-\epsilon}$, we have $\frac{\partial
J(\Theta_{d})}{\partial\theta_{i}}>0$ regardless of the value of $\Theta_{d}$.
Through \eqref{eq:param_update}, $\theta_{i}$ for every node $i$ will
eventually be reduced to the optimal value zero.
$\blacksquare$

\begin{remark} 
The result of Theorem 1 is consistent with, but more general
than, a similar result in \cite{yu2017optimal} where homogeneous targets are
assumed ($A_{i}=A$ and $B_{i}=B$ for all node $i$). The convergence of $\nabla
R(t)$ is related to the coefficients $A_{i}$ and $B_{i}$, $i=1,\ldots,M$. From
elementary queueing theory, a basic requirement for stability is $A_{i}<B_{i}$
for each node, which in turn implies the existence of $t_{K}$ in Assumption 1.
Moreover, if $T$ is sufficiently large and $\epsilon$ is arbitrarily small,
$\lim_{T\rightarrow\infty}\frac{\partial J(\Theta_{d})}{\partial\theta_{m}%
}\rightarrow1$ .
\end{remark}

\section{Simulation Examples}

\textbf{One agent, two targets}. We provide a simple one-agent example to
illustrate Theorem 1. Consider a controller with parameter vector $\Theta
_{d}=\left[  \theta_{1},\theta_{2}\right]  ^{\top}$. We track the evolution of
$\nabla R(t)=\left[  \frac{\partial R_{1}}{\partial\theta_{1}},\frac{\partial
R_{1}}{\partial\theta_{2}},\frac{\partial R_{2}}{\partial\theta_{1}}%
,\frac{\partial R_{2}}{\partial\theta_{2}}\right]  ^{\top}$ event by event.


\subsubsection{If the agent departs from target $1$}%

\[
\nabla R(\tau_{k}^{d_{1}^{+}}) = \hspace*{-2mm}%
\begin{bmatrix}
\frac{A_{1}}{A_{1} - B _{1}} & 0 & 0 & 0\\
0 & \frac{A_{1}}{A_{1} - B _{1}} & 0 & 0\\
0 & 0 & 1 & 0\\
0 & 0 & 0 & 1
\end{bmatrix}
\nabla R(\tau_{k}^{d_{1}^{-}}) +
\begin{bmatrix}
\frac{B_{1}}{B_{1} - A_{1}}\\
0\\
0\\
0
\end{bmatrix}
\]
For notational simplicity, denote the update matrix and update vector by
$\Lambda_{1}$ and $U_{1}$ respectively. We can write
\begin{equation}
\nabla R(\tau_{k}^{d_{1}^{+}}) = \Lambda_{1} \nabla R(\tau_{k}^{d_{1}^{-}}) +
U_{1} \label{eq:1A2T_case1}%
\end{equation}

\subsubsection{If the agent arrives at target $2$}%

\[
\nabla R(\tau_{k}^{a_{2}^{+}}) = \hspace*{-2mm}%
\begin{bmatrix}
1 & 0 & 0 & 0\\
0 & 1 & 0 & 0\\
\frac{B_{2}}{B_{1} - A_{1}} & 0 & 1 & 0\\
0 & \frac{B_{2}}{B_{1} - A_{1}} & 0 & 1
\end{bmatrix}
\nabla R(\tau_{k}^{a_{2}^{-}}) +
\begin{bmatrix}
0\\
0\\
\frac{B_{2}}{A_{1} - B_{1}}\\
0
\end{bmatrix}
\]
We denote this update by
\begin{equation}
\nabla R(\tau_{k}^{a_{2}^{+}}) = \Lambda_{2} \nabla R(\tau_{k}^{a_{2}^{-}}) +
U_{2} \label{eq:1A2T_case2}%
\end{equation}

\subsubsection{If the agent departs from target $2$}%

\[
\nabla R(\tau_{k}^{d_{2}^{+}}) =\hspace*{-2mm}
\begin{bmatrix}
1 & 0 & 0 & 0\\
0 & 1 & 0 & 0\\
0 & 0 & \frac{A_{2}}{A_{2} - B _{2}} & 0\\
0 & 0 & 0 & \frac{A_{2}}{A_{2} - B _{2}}%
\end{bmatrix}
\nabla R(\tau_{k}^{d_{2}^{-}}) +
\begin{bmatrix}
0\\
0\\
0\\
\frac{B_{2}}{B_{2} - A_{2}}%
\end{bmatrix}
\]
We denote this update by
\begin{equation}
\nabla R(\tau_{k}^{d_{2}^{+}}) = \Lambda_{3} \nabla R(\tau_{k}^{d_{2}^{-}}) +
U_{3} \label{eq:1A2T_case3}%
\end{equation}

\subsubsection{If the agent arrives at target $1$}%

\[
\nabla R(\tau_{k}^{a_{1}^{+}}) = \hspace*{-2mm}
\begin{bmatrix}
1 & 0 & \frac{B_{1}}{B_{2} - A_{2}} & 0\\
0 & 1 & 0 & \frac{B_{1}}{B_{2} - A_{2}}\\
0 & 0 & 1 & 0\\
0 & 0 & 0 & 1
\end{bmatrix}
\nabla R(\tau_{k}^{a_{1}^{-}}) +
\begin{bmatrix}
0\\
\frac{B_{1}}{A_{2} - B_{2}}\\
0\\
0
\end{bmatrix}
\]
We denote this update by
\begin{equation}
\nabla R(\tau_{k}^{a_{1}^{+}}) = \Lambda_{4} \nabla R(\tau_{k}^{a_{1}^{-}}) +
U_{4} \label{eq:1A2T_case4}%
\end{equation}

We initialize the agent at target $1$. The agent goes to target $2$ and back
to target $1$ so on so forth according to controller $\Theta_{d} = [\theta
_{1},\theta_{2}]^{\top}$. In each visiting cycle (from case 1 to case 4),
$\nabla R(t)$ is updated following the order from
\eqref{eq:1A2T_case1}-\eqref{eq:1A2T_case4}. We use $T_{k}$ to denote the beginning of the $k$-th cycle, and obtain
\begin{equation}%
\begin{split}
\nabla R(T_{k+1}^{-})\!= \!  &  \Lambda_{4}\! \left(  \Lambda_{3} \! \left(
\Lambda_{2} \!\left(  \Lambda_{1} \!\nabla R(T_{k}^{-}) + \! U_{1} \right)  +
U_{2} \right)  + U_{3} \right)  \!+\!U_{4}\\
=  &  \Lambda_{4}\Lambda_{3}\Lambda_{2}\Lambda_{1} \nabla R(T_{k}^{-}) +
\Lambda_{4}\Lambda_{3}\Lambda_{2} U_{1} + \Lambda_{4}\Lambda_{3} U_{2}\\
&  + \Lambda_{4} U_{3} + U_{4}\\
=  &  \Lambda\nabla R(T_{k}^{-}) + U
\end{split}
\label{eq:1A2T_sys}%
\end{equation}
where $\Lambda= \Lambda_{4}\Lambda_{3}\Lambda_{2}\Lambda_{1}$ and $U =
\Lambda_{4}\Lambda_{3}\Lambda_{2} U_{1} + \Lambda_{4}\Lambda_{3} U_{2} +
\Lambda_{4} U_{3} + U_{4} $ and the initial value $\nabla R(T_{0}) =
[0,0,0,0]^{\top}$.

\begin{figure}[h]
\centering
\includegraphics[width=0.9\linewidth]{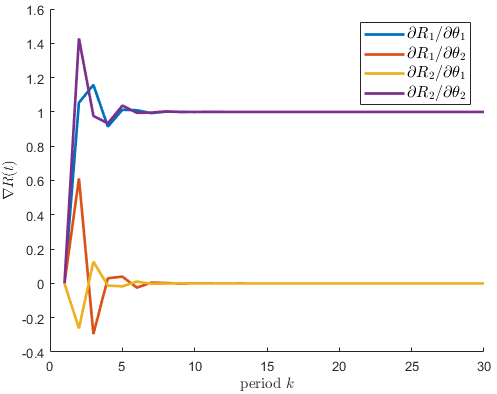}
\includegraphics[width=0.9\linewidth]{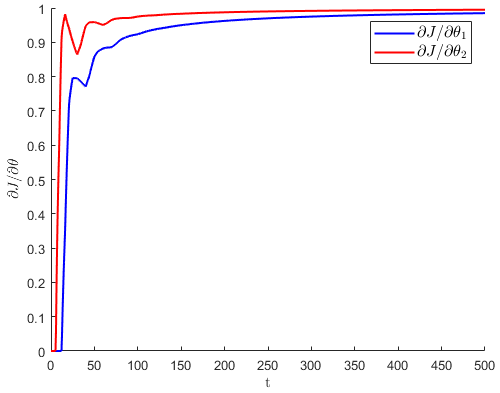}
\caption{{\protect\small Top plot: convergence of $\nabla R(t)$ to the
equilibrium $[1,0,0,1]^{\top}$. Bottom plot: convergence of $\partial
J/\partial\Theta_{d}$ to $[1,1]^{\top}$.}}%
\label{fig:Sim1_gradient_convergence}%
\end{figure}

Solving $\nabla R_{e} = (I- \Lambda)^{-1} U$, we obtain the only equilibrium $\nabla R_{e} = [1, 0, 0, 1]^{\top}$ of the system
\eqref{eq:1A2T_sys}. $\nabla R(t)$ converges to that
equilibrium asymptotically (see Fig. \ref{fig:Sim1_gradient_convergence}). The result matches with our analysis in the proof of
Theorem \ref{thm:thm1}. Moreover, convergence to the equilibrium simply
requires the eigenvalues of $\Lambda$ in \eqref{eq:1A2T_sys} lie within the
unit circle of the complex plane. 

\begin{figure}[h]
\centering
\includegraphics[width=0.95\linewidth]{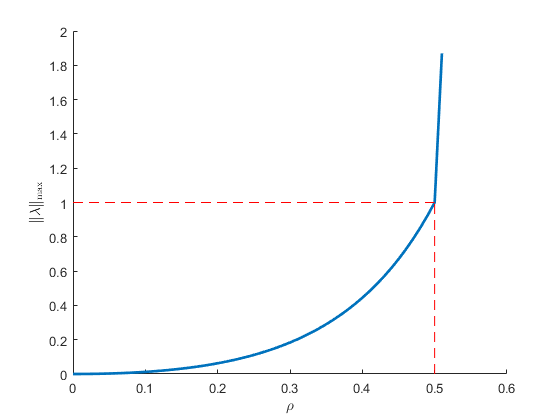}
\caption{{\protect\small Monotonic increasing of $\Vert\lambda\Vert_{\max}$ as the increase of $\rho$.}}%
\label{fig:eigenvalue_1A2T_h}%
\end{figure}

If the two targets are homogeneous ($A = A_{1} = A_{2}$ and $B = B_{1} =
B_{2}$), the convergence of $\nabla R(t)$ is only determined by the ratio
$\rho=A/B$. Using this ratio, we solve the eigenvalues of the system in
\eqref{eq:1A2T_sys} and obtain
\begin{align}
\bm\lambda=
\begin{bmatrix}
\frac{2\rho^{4} - 6\rho^{3}+7\rho^{2} -2\rho+ \sqrt{\rho^{3}(2\rho
-1)(2\rho^{2}-5\rho+ 4)}}{2(\rho- 1)^{4}}\\
\frac{2\rho^{4} - 6\rho^{3}+7\rho^{2} -2\rho+ \sqrt{\rho^{3}(2\rho
-1)(2\rho^{2}-5\rho+ 4)}}{2(\rho- 1)^{4}}\\
\frac{2\rho^{4} - 6\rho^{3}+7\rho^{2} -2\rho- \sqrt{\rho^{3}(2\rho
-1)(2\rho^{2}-5\rho+ 4)}}{2(\rho- 1)^{4}}\\
\frac{2\rho^{4} - 6\rho^{3}+7\rho^{2} -2\rho- \sqrt{\rho^{3}(2\rho
-1)(2\rho^{2}-5\rho+ 4)}}{2(\rho- 1)^{4}}%
\end{bmatrix}
\end{align}
Figure \ref{fig:eigenvalue_1A2T_h} shows the largest norm eigenvalue
$\Vert\lambda\Vert_{\max}$ increases monotonically as $\rho$ increases and
that $\Vert\lambda\Vert_{\max} =1$ at $\rho= 1/2$. The
convergence of $\nabla R(t)$ requires $\rho< 1/2$. 

Setting $\rho= 0.3$, we verify the convergence of both $\nabla
R(t)$ and $\partial J / \partial\Theta_{d}$ as shown in
Fig.\ref{fig:Sim1_gradient_convergence}. The results match with our analysis in Theorem 1. 

\textbf{Multi-agent cases: a counterexample to Theorem 1.} Theorem $1$ asserts
that an agent visiting a node should reduce the uncertainty state to zero
before moving to the next node. This property, however, does not apply to
multi-agent cases. This is not surprising because when two or more agents are
visiting a node, the allocation of agents to nodes may be improved if one
agent leaves the node before reducing its uncertainty state to zero and allow
other agents to complete this task.

\begin{figure}[h]
\centering
\includegraphics[width=0.8\linewidth]{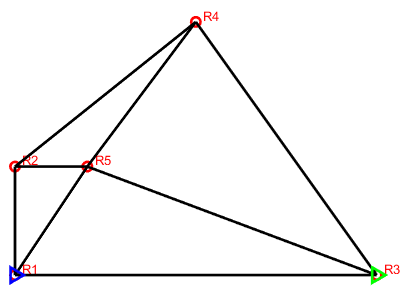}\caption{{\protect\small An
counter example with 2 homogeneous agents and 5 nodes to show $\theta
_{ii}^{a*} > 0$ for some agent $a$. }}%
\label{fig:counter_ex_2A5T}%
\end{figure}

Here we present a counterexample to Theorem 1 using two agents and
five nodes (see Fig. \ref{fig:counter_ex_2A5T}). Agents are initialized
at nodes 1 and 3 respectively and nodes are located at $X_{1}=(0,0)$, $X_{2}=(0,3)$, $X_{3}=(10,0)$, $X_{4}=(5,7)$, $X_{5}=(2,3)$ with uncertainty states $R_{i}(0)=0.5, A_i = 1, B_i = 10$ for $i=1,\ldots,5$. The initial thresholds are listed as follows:
\[
(\Theta^{1})^{0}=%
\begin{bmatrix}
16.34 & 5.31 & 5.18 & 1.74 & 0.72\\
2.87 & 1.02 & 18.56 & 22.13 & 24.55\\
23.76 & 9.93 & 9.80 & 23.82 & 8.49\\
12.05 & 5.83 & 4.56 & 23.28 & 17.67\\
21.81 & 21.04 & 18.59 & 10.39 & 9.05
\end{bmatrix}
\]

\[
(\Theta^{2})^{0}=
\begin{bmatrix}
0.88 & 22.13 & 3.33 & 10.81 & 22.28\\
21.38 & 22.60 & 17.45 & 0.45 & 22.96\\
16.43 & 0.26 & 9.96 & 17.29 & 1.83\\
19.14 & 1.86 & 22.08 & 11.74 & 1.14\\
13.85 & 6.12 & 4.53 & 3.21 & 10.96
\end{bmatrix}
\]

The final thresholds after convergence of (\ref{eq:param_update}), in this
case 300 iterations, are as follows:%

\[
(\Theta^{1})^{300}=
\begin{bmatrix}
0 & 2.38 & 5.18 & 1.74 & 0\\
5.70 & 0 & 18.56 & 23.76 & 23.94\\
23.76 & 9.93 & 7.25 & 23.82 & 8.51\\
12.05 & 5.83 & 4.56 & 8.06 & 17.67\\
21.70 & 21.04 & 20.09 & 16.07 & 0.23
\end{bmatrix}
\]

\[
(\Theta^{2})^{300}=
\begin{bmatrix}
0.88 & 22.13 & 3.33 & 10.81 & 22.24\\
21.37 & 19.27 & 17.45 & 0.17 & 25.24\\
16.44 & 1.15 & 0.02 & 15.23 & 2.21\\
19.14 & 1.86 & 22.08 & 0.01 & 0.87\\
13.85 & 6.12 & 2.27 & 2.21 & 0.00
\end{bmatrix}
\]

\begin{figure}[h]
\centering
\includegraphics[width=0.48\linewidth]{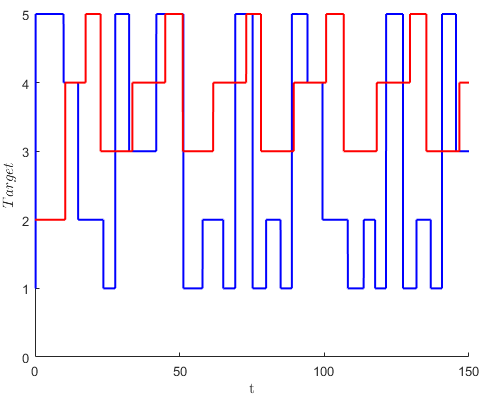}
\includegraphics[width=0.48\linewidth]{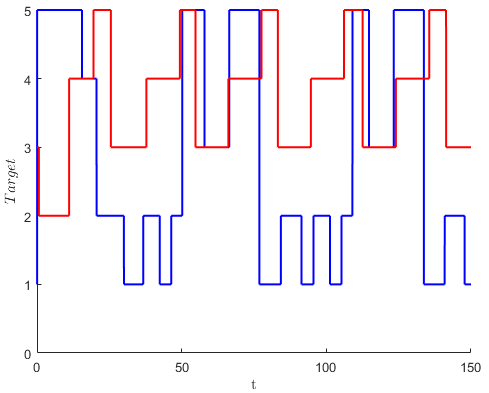}\caption{{\protect\small Left plot: the visiting sequence under the initial parameter. Right plot:
the sequence under the optimized parameters after 300 iterations of gradient
descent. In both plots, blue lines indicate the sequence of agent 1 and red
lines indicate the sequence of agent 2.}}%
\label{fig:tar_visit_seq_init_vs_opt}%
\end{figure}

\begin{figure}[h]
\centering
\includegraphics[width=1.0\linewidth]{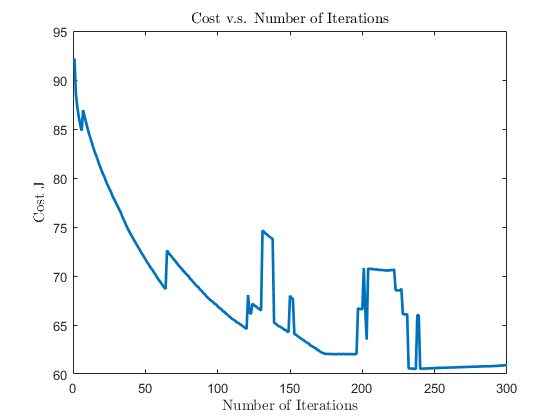}\caption{{\protect\small Cost
versus the number of iterations for the example of 2 agents and 5 nodes.}}%
\label{fig:Sim2_J_vs_iter}%
\end{figure}

The diagonal entries of the final parameter matrices for both agents are:
$\Theta_{d}^{1\ast}=[0,0,7.25,8.06,0.23]^{\top}$ and $\Theta_{d}^{2\ast
}=[0.88,19.27,0.02,0.01,0]^{\top}$ which do not satisfy the structure given in
Theorem $1$ as opposed to one-agent cases. In addition, the target visiting
sequences are adjusted on line during the optimization process. For instance,
the visiting sequence of agent $1$ is adjusted from initially being
$1-5-4-2-1-5-\ldots$ to $1-5-4-2-1-2-\ldots$ after 300 iterations as shown in
Fig. \ref{fig:tar_visit_seq_init_vs_opt}. Since agents may adjust their
visiting sequences asynchronously, the cost in the multi-agent cases may
fluctuate during the optimization process as shown in Fig.
\ref{fig:Sim2_J_vs_iter}

\textbf{Threshold-based policy versus dynamic programming. }We present a small
example to compare the performance of the threshold-based policy with a
classical dynamic programming solution of (\ref{eq:costfunction}) adapted to
the graph topology using value iteration.

\begin{figure}[h]
\centering
\includegraphics[width=1.0\linewidth]{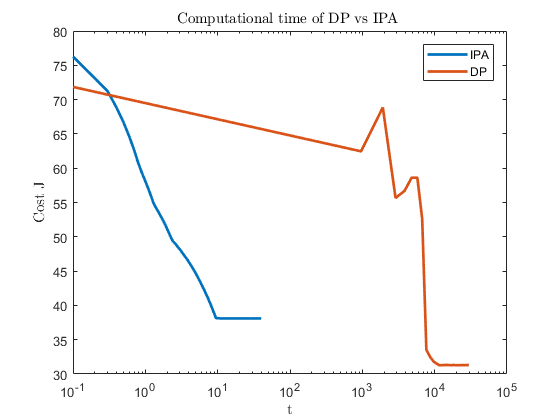}\caption{{\protect\small Cost
versus computational time (in log scale). The blue line shows the result of IPA
with the final cost $J^{\star}_{\text{IPA}} = 36.20$ and the orange line shows the result
of dynamic programming with the final cost $J^{\star}_{\text{DP}} = 31.15$ .}}%
\label{fig:sim_IPA_vs_DP}%
\end{figure}

A single agent is initialized at $(0,0)$ to persistently monitor four targets
located at $X_{1}=(0,0)$, $X_{2}=(4,0)$, $X_{3}=(4,4)$, $X_{4}=(0,4)$ (see
Fig. \ref{fig:ProbForm}) for $T=100$ seconds. The parameters in the
uncertainty dynamics \eqref{eq:DynR} are $A_{i}=1$, $B_{i}=20$, for
$i=1,\ldots,4$ and initial values are $R_{1}(0)=19,R_{2}(0)=14,R_{3}(0)=9$ and
$R_{4}(0)=4$. Using dynamic programming, the value function converges after
$15$ iterations and the final cost is $J_{\text{DP}}^{\star}=31.15$. However,
the number of states in the system consisting of 1 agent and 4 targets
$(\mathbf{s}(t),\mathbf{R}(t))$ is about $2.5\ast10^{9}$ discretized by
integers over 100 seconds. The running time is about 16 minutes per value
iteration using a computer with Intel(R) Core(TM) i7-7700 CPU @3.60GHZ
processor. Obviously, this method does not scale well in the number of states.
On the other hand, the solution obtained by optimizing the threshold-based
policy using the IPA approach is slightly higher, but the computational
complexity is reduced by several orders of magnitude as shown in Fig.
\ref{fig:sim_IPA_vs_DP}. After 300 iterations of gradient descent through
\eqref{eq:param_update} (about 30 seconds in total running time on the same
computer), the cost is reduced to $J_{\text{IPA}}^{\star}=36.20$.

\section{Conclusions}

The optimal multi-agent persistent monitoring problem involves the planning of
agent trajectories defined both by the sequence of nodes (targets) to be
visited and the amount of time spent by agents at each node. We have
considered a class of distributed parametric controllers through which the
agents control their visit sequence and dwell times at nodes using threshold
parameters associated with the node uncertainty states. We use Infinitesimal
Perturbation Analysis (IPA) to determine on line (locally) optimal threshold
parameters through gradient descent methods and thus obtain optimal
controllers within this family of threshold-based policies. In the one-agent
case we show the optimal strategy is for the agent to reduce the uncertainty
of a node to zero before moving to the next node. Compared with dynamic
programming solutions (in the limited instances when these are feasible), our
threshold-based parametric controller is effective and the computational
complexity is reduced by orders of magnitude. In future work, richer families
of threshold-based controllers can be developed by considering
multi-step-look-ahead policies and by identifying structural properties
therein which give us insight to the trade-off between exploitation and
exploration over multiple steps in persistent monitoring tasks.

\bibliographystyle{IEEEtran}
\bibliography{library}

\end{document}